\documentclass[reqno]{amsart}
\usepackage{graphicx}
\usepackage{color}
\usepackage{OddPatterns}
\pdfoutput=1

\let\TSp\thinspace
\def\DSp{\thinspace\thinspace}
\def\lbr{\raise 1pt\hbox{[}}
\def\rbr{\raise 1pt\hbox{]}}

\def\gobble#1{}
\def\hline{\multispan\numcolumns\hrulefill\cr}
\def\torframe#1{\vtop{\vbox{\hrule\hbox{\vrule\strut #1\vrule}}\hrule}}
\def\dblhline{\hline height 0.16667em\gobble&\emptyline\cr\hline}

\newbox\tablebox

\def\KhoHo{{\tt KhoHo}\space}
\def\tbb{\overline{tb}}

\def\hw{\protect\operatorname{hw}}
\def\ohw{\protect\operatorname{ohw}}
\def\thw{\widetilde{\protect\operatorname{hw}}}
\def\tohw{\widetilde{\protect\operatorname{ohw}}}
\let\lra\longrightarrow
\def\odd{\mathrm{odd}}
\def\oCalC{\CalC_{\odd}}
\def\oCalH{\CalH_{\odd}}
\def\od{d_{\odd}}

\theoremstyle{OVdefinition}
\newtheorem{example}[thm]{Example}

\begin{document}
\title{Patterns in odd Khovanov homology}
\author[A.~Shumakovitch]{Alexander Shumakovitch}
\address{Department of Mathematics, The George Washington University,
Phillips Hall, 801\ 22nd St. NW, Suite \#739, Washington, DC 20052, U.S.A.}
\email{Shurik@gwu.edu}
\thanks{The author was  partially supported by NSF grant DMS--0707526}
\begin{abstract}
We investigate properties of the odd Khovanov homology, compare and
contrast them with those of the original (even) Khovanov homology, and discuss
applications of the odd Khovanov homology to other areas of knot theory and
low-dimensional topology. We show that it provides an effective upper bound on
the Thurston-Bennequin number of Legendrian links and can be used to detect
quasi-alternating knots. A potential application to detecting transversely
non-simple knots is also mentioned.
\end{abstract}
\maketitle

\section{Introduction}
Khovanov homology is a special case of {\em categorification}, a novel
approach to construction of knot (or link) invariants that is being actively
developed over
the last decade after a seminal paper~\cite{Kh-Jones} by Mikhail Khovanov. The
idea of categorification is to replace a known polynomial knot (or link)
invariant with a family of chain complexes, such that the coefficients of the
original polynomial are the Euler characteristics of these complexes. Although
the chain complexes themselves depend heavily on a diagram that represents the
link, their homology depend solely on the isotopy class of the link. Khovanov
homology categorifies the Jones polynomial~\cite{Jones}.

More specifically, let $L$ be an oriented link in $\R^3$ represented by a
planar diagram $D$ and let $J_L(q)$ be a version of the Jones polynomial of
$L$ that satisfies the following identities (called the {\em Jones skein
relation} and {\em normalization}):
\begin{equation}\label{eq:Jones-skein}
-q^{-2}J_{\includegraphics[scale=0.45]{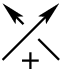}}(q)
+q^2J_{\includegraphics[scale=0.45]{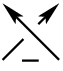}}(q)
=(q-1/q)J_{\includegraphics[scale=0.45]{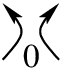}}(q);
\qquad
J_{\includegraphics[scale=0.45]{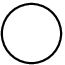}}(q)=q+1/q.
\end{equation}
The skein relation should be understood as relating the Jones polynomials of
three links whose planar diagrams are identical everywhere except in a small
disk, where they are different as depicted in~\eqref{eq:Jones-skein}. The
normalization fixes the value of the Jones polynomial on the trivial knot.
$J_L(q)$ is a Laurent polynomial in $q$ for every link $L$ and is
completely determined by its skein relation and normalization.
In our paper we also make use of another version of the Jones polynomial,
denoted $\tJ_L(q)$, that satisfies the same skein
relation~\eqref{eq:Jones-skein} but is normalized to equal~$1$ on the trivial
knot.

In~\cite{Kh-Jones} Mikhail Khovanov assigned to $D$ a family of Abelian groups
$\CalH^{i,j}(L)$ whose isomorphism classes depend on the isotopy class of $L$
only. These groups are defined as homology groups of an appropriate (graded)
chain complex $\CalC^{i,j}(D)$ with integer coefficients.
Groups $\CalH^{i,j}(L)$ are
nontrivial for finitely many values of the pair $(i,j)$ only. The gist of the
categorification is that the graded Euler characteristic of the Khovanov chain
complex equals $J_L(q)$:
\begin{equation}\label{eq:KhEuler-Jones}
J_L(q)=\sum_{i,j}(-1)^iq^jh^{i,j}(L),
\end{equation}
where $h^{i,j}(L)=\rk(\CalH^{i,j}(L))$, the Betti numbers of $\CalH$.
The reader is referred to~\cite{BN-first,Kh-Jones} for detailed treatment.

In 2007, Ozsv\'ath, Rasmussen and Szab\'o introduced~\cite{Khovanov-odd} an
{\em odd} version of the Khovanov homology. The odd Khovanov homology equals
the original (even) one modulo $2$ and, in particular, categorifies the same
Jones polynomial. On the other hand, the odd and even homology theories often
have drastically different properties (see Sections~\ref{sec:Kh-odd}
and~\ref{sec:Kh-comp} for details). The odd Khovanov homology appears to be
one of the connecting links between Khovanov and Heegaard-Floer homology
theories~\cite{OS-spectral}.

In this paper, we investigate properties of the odd Khovanov homology and
compare them to those of the even one. We also look at several applications
of the even and odd Khovanov homology theories, such as giving upper bounds on
the Thurston-Bennequin number of Legendrian links and detecting
quasi-alternating links, and observe that the odd homology often gives a
stronger result than the even one. Most of the experimental results that are
referred to in the paper were obtained with {\tt KhoHo}, a program by the
author for computing and studying Khovanov homology~\cite{Sh-KhoHo}.

This paper is organized as follows. In Section~\ref{sec:Kh-odd} we give a
definition of the odd Khovanov homology. Our exposition is self-contained, but
the reader is assumed to be familiar with the original Khovanov construction.
We compare even and odd Khovanov homology theories
with each other and list their basic properties in Section~\ref{sec:Kh-comp}. 
Section~\ref{sec:Kh-appl} is devoted to some of the more important
applications of the odd Khovanov homology to other areas of low-dimensional
topology.

The author is grateful to Mikhail Khovanov and Lenhard Ng for many fruitful
discussions during the work on this paper. He is also thankful to the referee
for several helpful comments and suggestions on the draft version of this
paper.

\section{Definition of the odd Khovanov homology}\label{sec:Kh-odd}
In this section we give a brief outline of the odd Khovanov homology theory
following~\cite{Khovanov-odd}. Our setting is slightly
more general than the one in the Introduction as we allow different
coefficient rings, not only $\Z$.

\subsection{Algebraic preliminaries}
Let $R$ be a commutative ring with unity. In this paper, we are mainly
interested in the cases when $R=\Z$, $\Q$, or $\Z_2$. If $M$ is a graded
$R$-module, we denote its {\em homogeneous component of degree $j$} by $M_j$.
For an integer $k$, the {\em shifted module} $M\{k\}$ is defined as having
homogeneous components $M\{k\}_j=M_{j-k}$. In the case when $M$ is free and
finite dimensional, we
define its {\em graded dimension} as the power series
$\dim_q(M)=\sum_{j\in\Z}q^j\dim(M_j)$ in variable $q$. Finally, 
if $(\CalC,d)=\cdots\lra\CalC^{i-1}\stackrel{d^{i{-}1}}{\lra}\CalC^i
\stackrel{d^i}{\lra}\CalC^{i+1}\lra\cdots$ is a (co)chain complex of graded
free $R$-modules such that all differentials $d^i$ are graded of degree $0$
with respect to the internal grading,
we define its {\sl graded Euler characteristic} as
$\chi_q(\CalC)=\sum_{i\in\Z}(-1)^i\dim_q(\CalC^i)$.

\begin{rem}
One can think of a (co)chain complex of graded $R$-modules as a {\em bigraded
$R$-module} where the homogeneous components are indexed by pairs of numbers
$(i,j)\in\Z^2$. Under this point of view, the differentials are graded of
{\em bidegree} $(1,0)$.
\end{rem}

\subsection{Odd Khovanov chain complex}\label{sec:Kh-complex} Let $L$ be an
oriented link and $D$ its planar diagram. We assign a number $\pm1$, called
{\em sign}, to every crossing of $D$ according to the rule depicted in
Figure~\ref{fig:crossing-signs}. The sum of these signs over all the crossings
of $D$ is called the {\em writhe number} of $D$ and is denoted by $w(D)$.

\begin{figure}
\captionindent 0.35\captionindent
\begin{minipage}{1.8in}
\centerline{\input{Xing_signs.pspdftex}}
\caption{Positive and negative crossings}
\label{fig:crossing-signs}
\end{minipage}
\hfill
\begin{minipage}{3.1in}
\centerline{\input{markers.pspdftex}}
\caption{Positive and negative markers and the corresponding resolutions of a
diagram.}
\label{fig:markers}
\end{minipage}
\end{figure}

Every crossing of $D$ can be {\em resolved} in two different ways according to
a choice of a {\em marker}, which can be either {\em positive} or {\em
negative}, at this crossing (see Figure~\ref{fig:markers}). A collection of
markers chosen at every crossing of a diagram $D$ is called a {\em (Kauffman)
state} of $D$. For a diagram with $n$ crossings, there are, obviously, $2^n$
different states. Denote by $\Gs(s)$ the the number of positive markers minus
the number of negative ones in a given state $s$. Define
\begin{equation}\label{eq:state-ij}
i(s)=\frac{w(D)-\Gs(s)}2,\qquad j(s)=\frac{3w(D)-\Gs(s)}2.
\end{equation}
Since both $w(D)$ and $\Gs(s)$ are congruent to $n$ modulo $2$, $i(s)$ and
$j(s)$ are always integers. For a given state $s$, the result of
the resolution of $D$ at each crossing according to $s$ is a family $D_s$ of
disjointly embedded circles. Denote the number of these circles by $|D_s|$.

For each state $s$ of $D$, we assign a free graded $R$-module $\GL(s)$ as
follows. Label all circles from the resolution $D_s$ by some independent
variables, say,
$X^s_1,X^s_2,\,\dots\,,X^s_{|D_s|}$ and let
$V_s=V(X^s_1,X^s_2,\,\dots\,,X^s_{|D_s|})$ be
a free $R$-module generated by them. We define $\GL(s)=\GL\!^*(V_s)$, the
exterior algebra of $V_s$. Then $\GL(s)=\GL\!^0(V_s)\oplus\GL\!^1(V_s)\oplus
\cdots\oplus\GL\!^{|D_s|}(V_s)$ and we grade $\GL(s)$ by specifying
$\GL(s)_{|D_s|-2k}=\GL\!^k(V_s)$ for each $0\le k\le|D_s|$, where
$\GL(s)_{|D_s|-2k}$ is the homogeneous component of $\GL(s)$ of degree
$|D_s|-2k$. It is an easy exercise for the reader to check that
$\dim_q(\GL(s))=(q+q^{-1})^{|D_s|}$.

Let $\oCalC^i(D)=\bigoplus_{i(s)=i}\GL(s)\{j(s)\}$ for each $i\in\Z$. In
order to arrange these modules into a graded chain complex $\oCalC(D)$,
we need to define a (graded) differential $\od^i:\oCalC^i(D)\to\oCalC^{i+1}(D)$
of degree $0$. But even before the chain complex structure on $\oCalC(D)$ is
defined, its (graded) Euler characteristic makes sense. Similarly to the case
of the even Khovanov homology, one can easily verify that
$\chi_q(\oCalC(D))=J_L(q)$. In fact, $\oCalC(D)\simeq\CalC(D)$ as bigraded
$R$-modules~\cite{Khovanov-odd}.

\begin{figure}
\centerline{\input{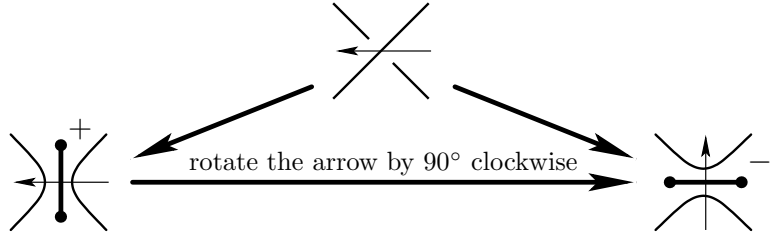}}
\caption{Choice of arrows at the diagram crossings}\label{fig:odd-arrows}
\end{figure}

To define a differential on $\oCalC$, we introduce an additional structure,
a choice of an arrow at each crossing of $D$ that is parallel to the {\em
negative} marker at that crossing (see Figure~\ref{fig:odd-arrows}). There
are obviously $2^n$ such choices. For every state $s$ on $D$, we place arrows
that connect two branches of $D_s$ near each (former) crossing according to
the rule from Figure~\ref{fig:odd-arrows}.

\begin{figure}
\centerline{\vbox{\halign{#\hfill\cr
\hskip-0.5em\input{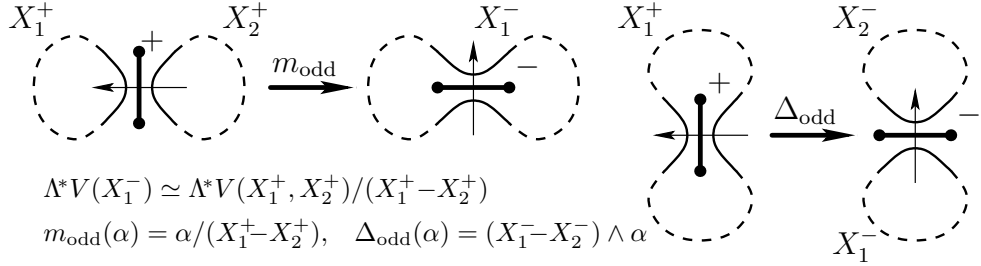}\cr
\noalign{\vskip -1.2\baselineskip}
\hbox{\vbox to0pt{\vss\halign{#\hfill\cr
$\GL\!^*V(X^-_1)\simeq\GL\!^*V(X^+_1,X^+_2)/(X^+_1{-}X^+_2)$\cr
\noalign{\medskip}
$m_{\odd}(\Ga)=\Ga/(X^+_1\!\!{-}X^+_2)$,\quad
$\GD_{\odd}(\Ga)=(X^-_1\!\!{-}X^-_2)\wedge\Ga$\cr
\noalign{\medskip}}}}\cr
}}}
\caption{Diagram resolutions corresponding to adjacent states and
maps between the free $R$-modules assigned to these resolutions}
\label{fig:odd-diff}
\end{figure}

Let $s_+$ and $s_-$ be two states of $D$ that differ at a single
crossing, where $s_+$ has a positive marker while $s_-$ has a negative one.
We call two such states {\em adjacent}. In this case, $\Gs(s_-)=\Gs(s_+)-2$
and, consequently, $i(s_-)=i(s_+)+1$ and $j(s_-)=j(s_+)+1$. Consider now the
resolutions of $D$ corresponding to $s_+$ and $s_-$. One can
readily see that $D_{s_-}$ is obtained from $D_{s_+}$ by either merging two
circles into one or splitting one circle into two (see
Figure~\ref{fig:odd-diff}). All the circles that do not pass through the
crossing at which $s_+$ and $s_-$ differ, remain unchanged.
We define $d_{s_+:s_-}:\GL(s_+)\to\GL(s_-)$ as either $m_{\odd}$ or
$\GD_{\odd}$ depending on whether the circles merge or split, where $m_{\odd}$
and $\GD_{\odd}$ are as follows.

If $s_-$ is obtained from $s_+$ by merging two circles together, then
we have that $\GL(s_-)\simeq\GL(s_+)/(X^+_1-X^+_2)$,
where $X^+_1$ and $X^+_2$ are the generators
of $V_{s_+}$ corresponding to the two merging circles, as depicted in
Figure~\ref{fig:odd-diff}. We define $m_{\odd}:\GL(s_+)\to\GL(s_-)$ to be
this isomorphism composed with the projection
$\GL(s_+)\to\GL(s_+)/(X^+_1-X^+_2)$.

The case when one circle splits into two is more interesting. Let $X^-_1$ and
$X^-_2$ be the generators of $V_{s_-}$ corresponding to these two circles such
that the arrow points from $X^-_1$ to $X^-_2$ (see Figure~\ref{fig:odd-diff}).
Now for each generator $X^+_k$ of $V_{s_+}$, we define
$\GD_{\odd}(X^+_k)=(X^-_1-X^-_2)\wedge X^-_{\eta(k)}$ where $\eta$ is the
correspondence between circles in $D_{s_+}$ and $D_{s_-}$. While $\eta(1)$ can
equal either $1$ or $2$, this choice does not affect $\GD_{\odd}(X^+_1)$ since
$(X^-_1-X^-_2)\wedge X^-_2=X^-_1\wedge X^-_2=-X^-_2\wedge X^-_1=
(X^-_1-X^-_2)\wedge X^-_1$.

We need one more ingredient in order to finish the definition of the
differential on $\oCalC(D)$. Namely, to each adjacent pair of states
$(s_+,s_-)$, we assign a sign $\Ge(s_+,s_-)\in\{\pm1\}$. This choice of a sign
for each adjacent pair is called an {\em edge assignment}. The etymology of
this term will become clear in a few paragraphs. Finally, let
$\od^i=\sum_{(s_+,s_-)}\Ge(s_+,s_-)d_{s_+:s_-}$, where $(s_+,s_-)$ runs over
all adjacent pairs of states with $i(s_+)=i$. Our goal is to choose an edge
assignment in such a way that $\od:\CalC(D)\to\CalC(D)$ is indeed a
differential, that is, $\od^{i+1}\circ \od^i=0$ for every $i$. The following
Theorem guarantees us that it is always possible.

\begin{thm}[Ozsv\'ath--Rasmussen--Szab\'o~\cite{Khovanov-odd}]
\label{thm:ORS-odd}
It is possible to assign a sign $\Ge(s_+,s_-)$ to each adjacent pair of states
$(s_+,s_-)$ in such a way that $\oCalC(D;R)$ equipped with the differential
$\od$ defined above becomes a graded (co)chain complex.
The homology $\oCalH(L;R)$ of $\oCalC(D;R)$ does not depend on the
choice of arrows at the crossings, the choice of edge assignment,
and some other
choices needed in the construction. Moreover, the isomorphism class of
$\oCalH(L;R)$ is a link invariant that categorifies $J_L(q)$, a version of the
Jones polynomial defined by~\eqref{eq:Jones-skein}.
\end{thm}

\begin{defin}[\cite{Khovanov-odd}]
$\oCalH(L;R)$ is called the {\em odd Khovanov homology} of $L$ over a ring of
coefficients $R$. If $R$ is omitted from the notation, integer coefficients
are assumed.
\end{defin}

\begin{attn}\label{attn:odd-Z2}
By comparing the definitions of $\oCalC(D;\Z_2)$ and $\CalC(D;\Z_2)$, it
is easy to see that they are isomorphic as graded chain complexes (since the
signs do not matter modulo $2$). It follows that
$\oCalH(D;\Z_2)\simeq\CalH(D;\Z_2)$ as well.
\end{attn}

\begin{figure}
\centerline{\input{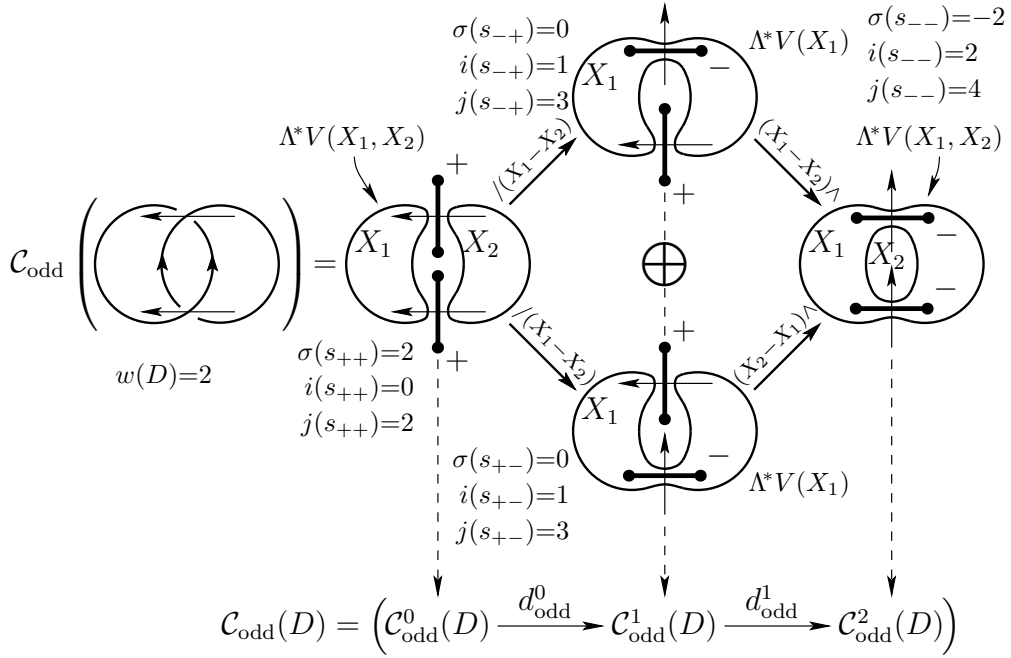}}
\caption{Odd Khovanov chain complex for the Hopf link}\label{fig:Kh-Hopf-odd}
\end{figure}

\begin{example}\label{ex:Hopf-link}
Figure~\ref{fig:Kh-Hopf-odd} shows the odd Khovanov chain complex for the Hopf
link with the indicated orientation. The diagram has two positive crossings,
so its writhe number is $2$. Let $s_{\pm\pm}$ be the four possible resolutions
of this diagram, where each ``$+$'' or ``$-$'' describes the sign of the
marker at the corresponding crossing. By looking at
Figure~\ref{fig:Kh-Hopf-odd}, one easily determines that
\begin{align*}
\oCalC^0(D)&=\GL(s_{++})\{2\}=\GL\!^*V(X_1,X_2)\{2\};\\
\oCalC^1(D)&=\GL(s_{+-})\{3\}\oplus \GL(s_{-+})\{3\}=
(\GL\!^*V(X_1)\oplus \GL\!^*V(X_1))\{3\};\\
\oCalC^2(D)&=\GL(s_{--})\{4\}=\GL\!^*V(X_1,X_2)\{4\}.
\end{align*}

It is convenient to arrange the four resolutions in the corners of a square
placed in the plane in such a way that its diagonal from $s_{++}$ to $s_{--}$
is horizontal. Then the edges of this square correspond to the maps between
the adjacent states (see Figure~\ref{fig:Kh-Hopf-odd}).
\end{example}

In general, $2^n$ resolutions of a diagram $D$ with $n$ crossings can be
arranged into an $n$-dimensional {\em cube of resolutions}, where vertices
correspond to the $2^n$ states of~$D$. The edges of this cube connect adjacent
pairs of states and can be oriented from $s_+$ to $s_-$. Every edge is assigned
either $m_{\odd}$ or $\GD_{\odd}$ with the sign $\Ge(s_+,s_-)$, as described
above. The differential $\od^i$ restricted to each summand $\GL(s)$
with $i(s)=i$ equals the sum of all the maps assigned to the edges that
originate at~$s$.

\begin{rem}
If one forgets about the signs on the edges, then each square (that is, a
$2$-dimensional face) in the cube of resolutions is either commutative, or
anti-commutative, or both. The latter case means that both double-composites
corresponding to the square are trivial. This is a major departure from the
situation in the even case, where each square is commutative~\cite{Kh-Jones}.
The role of the signs $\Ge(s_+,s_-)$ is to make all of the squares
anti-commutative.
\end{rem}

\begin{rem}
In the example~\ref{ex:Hopf-link}, the only square of resolutions is 
anti-commutative, so no adjustment of signs is needed.
\end{rem}

\begin{rem}
There is no explicit construction for finding edge assignments for
the cube of resolutions in the case of the odd Khovanov chain complex. 
Theorem~\ref{thm:ORS-odd} only ensures that signs exist. On the other hand,
such an assignment is fairly straightforward for the even Khovanov chain
complex.
\end{rem}

\subsection{Reduced odd Khovanov homology}

Similar to the even case~\cite{Kh-patterns}, there is a reduced version of the
odd Khovanov homology:
\begin{thm}[\cite{Khovanov-odd}]\label{thm:odd-reduced}
For every link $L$, there exists a bigraded abelian group $\tCalH_{\odd}(L;R)$
such that
$\oCalH(L;R)\simeq\tCalH_{\odd}(L;R)\{1\}\oplus \tCalH_{\odd}(L;R)\{-1\}$.
\end{thm}

$\tCalH_{\odd}(L;R)$ is called the {\em reduced odd} Khovanov homology of a
link $L$ over $R$. Theorem~\ref{thm:odd-reduced} implies that the reduced and
non-reduced odd homology determine each other completely. It is therefore
enough to consider the reduced version of the odd Khovanov homology only. This
contrasts with the even case where there are several examples of pairs of
knots that have the same rational Khovanov homology, but different rational
reduced ones.

\begin{rem}
One of the main differences between the even and odd reduced Khovanov
homologies is that the former depends on the choice of a link component, while
the latter does not. This difference disappears over $\Z_2$.
\end{rem}

\begin{thm}[\cite{Sh-torsion}]\label{thm:me-Z2}
$\CalH(L;\Z_2)\simeq\tCalH(L;\Z_2)\otimes_{\Z_2}A_{\Z_2}$ for every link $L$, 
where $A_{\Z_2}=\Z_2[X]/X^2$ is the algebra of truncated polynomials graded by
$\deg(1)=1$ and $\deg(X)=-1$. In particular, $\tCalH(L;\Z_2)$ does not depend
on the choice of a component of $L$.
\end{thm}

The odd Khovanov homology should provide an insight into interrelations
between Khovanov and Heegaard-Floer~\cite{OS-HF} homology theories. Its
definition was motivated by the following result.

\begin{thm}[Ozsv\'ath--Szab\'o~\cite{OS-spectral}]\label{thm:OS-spec}
For each link $L$ with a diagram $D$, there exists a spectral sequence with
$E^1=\tCalC(D;\Z_2)$ and $E^2=\tCalH(L;\Z_2)$ that converges to the
$\Z_2$-Heegaard-Floer homology $\widehat{HF}(\GS(L);\Z_2)$ of the double
branched cover $\GS(L)$ of $S^3$ along $L$.
\end{thm}

\begin{conjec}\label{conj:OS-spec}
There exists a spectral sequence that starts with $\tCalC_{\odd}(D;\Z)$
and $\tCalH_{\odd}(L;\Z)$ and converges to $\widehat{HF}(\GS(L);\Z)$.
\end{conjec}

\section{Comparison between even and odd Khovanov homology}
\label{sec:Kh-comp}

In this section we summarize the main differences and similarities in
properties exhibited by the even and odd versions of Khovanov homology and
list related constructions.

\setcounter{thm}{0}
\begin{attn}
Let $L$ be an oriented link and $D$ its planar diagram. Then
\begin{itemize}
\item $\oCalH(L;\Z_2)\simeq\CalH(L;\Z_2)$ and
$\tCalH_{\odd}(L;\Z_2)\simeq\tCalH(L;\Z_2)$.
\item $\chi_q(\CalH(L;R))=\chi_q(\oCalH(L;R))=J_L(q)$ and\\
$\chi_q(\tCalH(L;R))=\chi_q(\tCalH_{\odd}(L;R))=\tJ_L(q)$.
\item
$\CalH(L;\Z_2)\simeq\tCalH(L;\Z_2)\otimes_{\Z_2}\!A_{\Z_2}$~\cite{Sh-torsion}
and $\oCalH(L;R)\simeq\tCalH_{\odd}(L;R)\{1\}\oplus
\tCalH_{\odd}(L;R)\{-1\}$~\cite{Khovanov-odd}. On the other hand,
$\CalH(L;\Z)$ and $\CalH(L;\Q)$ do not split in general.
\item For links, $\tCalH(L;\Z_2)$ and $\tCalH_{\odd}(L;R)$ do not depend on the
choice of a component on which a base point lies. This is, in general, not the
case for $\tCalH(L;\Z)$ and $\tCalH(L;\Q)$.
\item If $L$ is a non-split alternating link, then $\CalH(L;\Q)$,
$\tCalH(L;R)$, and $\tCalH_{\odd}(L;R)$ are completely determined by the
Jones polynomial and signature of~$L$~\cite{Kh-patterns,Lee,Khovanov-odd}.
\item $\CalH(L;\Z_2)$ and $\oCalH(L;\Z)$ are invariant under
link mutations~\cite{Bloom-mutation,Wehrli-mutation2}. On the other hand,
$\CalH(L;\Z)$ is known not to be preserved under a mutation that
exchanges components of a link~\cite{Wehrli-mutation1} and under a
cabled (or genus~$2$) mutation~\cite{DGShT}. It is unclear whether
$\CalH(L;\Z)$ is invariant under a component-preserving mutation.
\item $\CalH(L;\Z)$ almost always has torsion (except for several special
cases), but usually of order $2$. The first knot with $4$-torsion is the
$(4,5)$-torus knot that has 15 crossings. The first known knot with
$3$-torsion is the $(5,6)$-torus knot with 24 crossings. On the other hand,
$\oCalH(L;\Z)$ was observed to have torsion of various orders even for
knots with relatively few crossings (see remark on
page~\pageref{rem:odd-torsion}), although orders $2$ and $3$ are the most
popular.
\item $\tCalH(L;\Z)$ has very little torsion. The first knot with torsion has
13 crossings. On the other hand, $\tCalH_{\odd}(L;\Z)$ has as much torsion as
$\oCalH(L;\Z)$ (see Theorem~\ref{thm:odd-reduced}).
\end{itemize}
\end{attn}

\begin{rem}
The properties above show that $\tCalH_{\odd}(L;\Z)$ behaves
similarly to $\tCalH(L;\Z_2)$ but not to $\tCalH(L;\Z)$. This is one of its
main features (see~\ref{thm:OS-spec} and~\ref{conj:OS-spec}).
\end{rem}

\subsection{Distinguishing knots}
The odd Khovanov homology appears to be somewhat stronger of a knot invariant
than the even one, although independent. Among $313230$ prime knots with up to
$15$ crossings (111528 alternating and 201702 non-alternating),
there are 4377 pairs of knots (8754 if counted with mirror images,
as none of these knots are amphicheiral)
that have the same even Khovanov homology but different odd ones. The first
knots that are distinguished by the odd homology but not the even one are
$11^n_{1}$ and $\overline{12}^n_{577}$. On the other hand, there are only 107
(214 if counted with mirror images) pairs of prime knots with up to $15$
crossings that have the same odd Khovanov homology but different even ones. The
first knots distinguished by the even homology but not the odd one are
$13^n_{2640}$ and $15^n_{124915}$. We don't know any special properties of
these knots.

\begin{rem}
\label{rem:knots-enum}
Throughout this paper we use the following notation for knots: knots with 10
crossings or less are numbered according to Rolfsen's knot
table~\cite{Rolfsen-book} and knots with 11 crossings or more are numbered
according to the knot table from Knotscape~\cite{Knotscape}. Mirror images of
knots from either table are denoted with a bar on top. For example,
$\overline{9}_{46}$ is the mirror image of the knot number 46 with 9 crossings
from Rolfsen's table and $16^n_{197566}$ is the non-alternating
knot number 197566 with 16 crossings from Knotscape's table.
\end{rem}

It is interesting to notice that the odd Khovanov homology is the same for all
{\em almost mutant knots}~\cite{DGShT} with up to $16$ crossings. There is
currently no explanation for this fact. Even homology, on the other hand, is
known to distinguish some of them.

\subsection{Homological thickness}
\begin{defin}\label{def:hw}
Let $L$ be a link. The {\em homological width} of $L$ over a ring $R$ is the
minimal number of adjacent diagonals $j-2i=const$ such that $\CalH(L;R)$ is
zero outside of these diagonals.
It is denoted by $\hw_R(L)$. The {\em reduced homological width}, $\thw_R(L)$
of $L$, {\em odd homological width}, $\ohw_R(L)$ of $L$, and {\em reduced odd
homological width}, $\tohw_R(L)$ of $L$ are defined similarly.
\end{defin}

\begin{attn}\label{attn:hom-width}
It follows from~\ref{thm:odd-reduced} that $\tohw_R(L)=\ohw_R(L)-1$. The same
holds true in the case of the even Khovanov homology over $\Q$:
$\thw_\Q(L)=\hw_\Q(L)-1$ (see~\cite{Kh-patterns}).
\end{attn}

\begin{defin}\label{def:thick-thin}
A link $L$ is said to be {\em homologically thin} over a ring $R$, or simply
$R$H-thin, if $\hw_R(L)=2$. $L$ is {\em homologically thick}, or $R$H-thick,
otherwise. We define {\em odd-homologically} thin and thick, or simply
$R$OH-thin and $R$OH-thick, links similarly. 
\end{defin}

\begin{thm}[Lee, Ozsv\'ath--Rasmussen--Szab\'o,
Manolescu--Ozsv\'ath~\cite{Lee,Khovanov-odd,QA-links}]\label{thm:alt-thin}
Quasi-alternating links (see Section~\ref{sec:quasi-alt} for the definition)
are $R$H-thin and $R$OH-thin for every ring $R$. In particular, this is
true for non-split alternating links.
\end{thm}

\begin{thm}[Khovanov~\cite{Kh-patterns}]\label{thm:adequate-thick}
Adequate links are $R$H-thick for every $R$.
\end{thm}

\begin{figure}
\centerline{\input knot-15n_41127-table-Red}
\par\medskip
$\hw_{\Q}=3,\qquad \thw_{\Q}=2,\qquad \hw_{\Z}=4,\qquad \thw_{\Z}=3$
\par\bigskip\bigskip
\centerline{\input knot-15n_41127-table-Odd}
\par\medskip
$\ohw_{\Q}=4,\qquad \tohw_{\Q}=3,\qquad \ohw_{\Z}=4,\qquad \tohw_{\Z}=3$

\caption{Integral reduced even Khovanov homology (above) and odd Khovanov
Homology (below) of the knot $15^n_{41127}$}\label{fig:odd-thick}
\end{figure}

\begin{figure}
\centerline{\resizebox{\hsize}{!}{\input knot-16n_197566-table}}
\par\medskip
The free part of $\CalH(16^n_{197566};\Z)$ is supported on diagonals 
$j-2i=7$ and $j-2i=9$. On the other hand, there is $2$-torsion on the diagonal
$j-2i=5$. Therefore, $16^n_{197566}$ is $\Q$H-thin, but $\Z$H-thick and
$\Z_2$H-thick.
\par\bigskip\bigskip\bigskip
\centerline{\resizebox{\hsize}{!}{\input knot-16n_-197566-table}}
\par\medskip
$\CalH(\overline{16}^n_{197566};\Z)$ is supported on diagonals 
$j-2i=-7$ and $j-2i=-9$.\\
But there is $2$-torsion on the diagonal $j-2i=-7$.\\
Therefore, $\overline{16}^n_{197566}$ is $\Q$H-thin and
$\Z$H-thin, but $\Z_2$H-thick.
\par\bigskip
\caption{Integral Khovanov homology of the knots $16^n_{197566}$ and
$\overline{16}^n_{197566}$}\label{fig:diff-thickness}
\end{figure}

\begin{attn}
Homological thickness of a link $L$ often does not depend on the base ring.
The first prime knot with $\hw_{\Q}(L)<\hw_{\Z_2}(L)$ and
$\hw_{\Q}(L)<\hw_{\Z}(L)$ is $15^n_{41127}$ with 15 crossings (see
Figure~\ref{fig:odd-thick}). The first prime
knot that is $\Q$H-thin but $\Z$H-thick, $16^n_{197566}$, has 16 crossings
(see Figure~\ref{fig:diff-thickness}). Its mirror image,
$\overline{16}^n_{197566}$ is both $\Q$H- and $\Z$H-thin but is $\Z_2$H-thick
with $\CalH^{-8,-21}(\overline{16}^n_{197566};\Z_2)\simeq\Z_2$, for example,
because of the Universal Coefficient Theorem. Also observe that
$\CalH^{9,25}(16^n_{197566};\Z)$ and
$\CalH^{-8,-25}(\overline{16}^n_{197566};\Z)$ have $4$-torsion, shown in a
small box in the tables.
\end{attn}

\begin{rem}
Tables in Figures~\ref{fig:odd-thick} and~\ref{fig:diff-thickness}
tabulate non-zero homology groups of the corresponding knots where the
$i$-grading is represented horizontally and the $j$-grading vertically. A
table entry of $\mathbf1$ or $\mathbf{1_2}$ means that the corresponding group
is $\Z$ or $\Z_2$, respectively. In general, an entry of the form
$\mathbf{a,b_c}$ would correspond to the group $\Z^a\oplus\Z_c^b$.
\end{rem}

\begin{attn}
Odd Khovanov homology is often thicker over $\Z$ than the even one. This
is crucial for applications (see Section~\ref{sec:Kh-appl}).
On the other hand, $\tohw_{\Q}(L)\le\thw_{\Q}(L)$ for all but one prime knot
with at most $15$ crossings. The homology for this knot, $15^n_{41127}$, is
shown in Figure~\ref{fig:odd-thick}. Please observe that
$\tCalH_{\odd}(15^n_{41127})$ has $3$-torsion (in gradings $(-2,-2)$ and
$(-1,0)$), while $\tCalH(15^n_{41127})$ has none.
\end{attn}

\def\Cpic#1{%
  \oCalC\!\left(\vcenter{\hbox{\includegraphics[scale=0.35]{#1}}}\right)\!}
\def\Hpic#1{%
  \oCalH\!\left(\vcenter{\hbox{\includegraphics[scale=0.35]{#1}}}\right)\!}
\subsection{Long exact sequence of the Khovanov homology}
One of the most useful tools in studying Khovanov homology is the long exact
sequence that categorifies Kauffman's {\em unoriented} skein relation for
the Jones polynomial~\cite{Kh-Jones}. A similar long exact sequence exists for
the odd version of the homology as well~\cite{Khovanov-odd}.
If we forget about the grading, then it
is clear from the construction from Section~\ref{sec:Kh-complex} that
$\,\Cpic{hsmooth_Xing-black}\,$ is a subcomplex of
$\,\Cpic{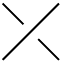}\,$ and $\,\Cpic{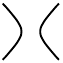}\simeq
\Cpic{just_Xing-black}/\Cpic{hsmooth_Xing-black}\,$ (see also
Figure~\ref{fig:Kh-Hopf-odd}). Here, 
$\vcenter{\hbox{\includegraphics[scale=0.35]{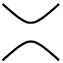}}}$ and
$\vcenter{\hbox{\includegraphics[scale=0.35]{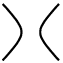}}}$ depict
link diagrams where a single crossing
$\vcenter{\hbox{\includegraphics[scale=0.35]{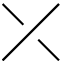}}}$ is
resolved in a negative or, respectively, positive direction. This results in a
short exact sequence of non-graded chain complexes:
\begin{equation}\label{eq:short-exact}
0\lra\Cpic{hsmooth_Xing-black}\stackrel{in}{\lra}\Cpic{just_Xing-black}
\stackrel{p}{\lra}\Cpic{vsmooth_Xing-black}\lra0,
\end{equation}
where $in$ is the inclusion and $p$ is the projection.

In order to introduce grading into~\eqref{eq:short-exact}, we need to consider
the cases when the crossing to be resolved is either positive or negative.
We get (see~\cite{Jake-comparison}):
\begin{equation}\label{eq:short-exact-orient}\begin{split}
0\lra\Cpic{hsmooth_Xing-black}\{2{+}3\Go\}[1{+}\Go]\stackrel{in}{\lra}
\Cpic{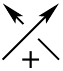}\stackrel{p}{\lra}
\Cpic{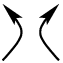}\{1\}\lra0,\\[\smallskipamount]
0\lra\Cpic{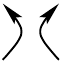}\{-1\}\stackrel{in}{\lra}
\Cpic{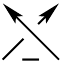}\stackrel{p}{\lra}
\Cpic{hsmooth_Xing-black}\{1{+}3\Go\}[\Go]\lra0,
\end{split}\end{equation}
where $\Go$ is the difference between the numbers of negative crossings in the
unoriented resolution
$\vcenter{\hbox{\includegraphics[scale=0.35]{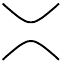}}}$ (it has to
be oriented somehow in order to define its Khovanov chain complex) and 
in the original diagram. The notation $\oCalC[k]$ is used to represent a
shift in the homological grading of a complex $\oCalC$ by $k$. The graded
versions of $in$ and $p$ are both homogeneous, that is, have bidegree $(0,0)$.

By passing to homology in~\eqref{eq:short-exact-orient}, we get the following
result.

\begin{thm}[Ozsv\'ath-Rasmussen-Szab\'o~\cite{Khovanov-odd}, see
also~\cite{Jake-comparison}]\label{thm:Kh-exact} The Khovanov homology is
subject to the following long exact sequences:
\begin{equation}\label{eq:Kh-exact}\begin{split}
\hskip -13pt
\cdots\lra\Hpic{smooth_Xing-nomark}\{1\}\stackrel{\partial}{\lra}
\Hpic{hsmooth_Xing-black}\{2{+}3\Go\}[1{+}\Go]\stackrel{in_*}{\lra}
\Hpic{pos_Xing-black}\stackrel{p_*}{\lra}
\Hpic{smooth_Xing-nomark}\{1\}\lra\cdots\\[\smallskipamount]
\hskip -13pt
\cdots\lra\Hpic{smooth_Xing-nomark}\{-1\}\stackrel{in_*}{\lra}
\Hpic{neg_Xing-black}\stackrel{p_*}{\lra}
\Hpic{hsmooth_Xing-black}\{1{+}3\Go\}[\Go]\stackrel{\partial}{\lra}
\Hpic{smooth_Xing-nomark}\{-1\}\lra\cdots
\end{split}\end{equation}
where $in_*$ and $p_*$ are homogeneous and $\partial$ is the connecting
differential and has bidegree $(1,0)$.
\end{thm}

\begin{rem}
Long exact sequences~\eqref{eq:Kh-exact} work equally well over any ring $R$
and for every version of Khovanov homology. This is both a blessing and a
curse. On one hand, this means that all of the properties of the even Khovanov
homology that are proved using these long exact sequences (and most of them
are) hold automatically true for the odd Khovanov homology as well. On the
other hand, this makes it very hard to find explanations to many differences
between these homology theories.
\end{rem}

\section{Applications of the odd Khovanov homology}\label{sec:Kh-appl}

In this section we collect some of the more prominent applications of
odd Khovanov homology. We make a special effort to compare its performance 
to even homology. Unfortunately, there appears to be no odd
version of the Rasmussen invariant~\cite{Jake-Milnor}. In fact, a knot might
not have any rational homology in the homological grading $0$ of odd
Khovanov homology at all, see Figure~\ref{fig:12n_475-hom}. Correspondingly,
odd homology cannot be used to provide bounds on the slice genus of a knot
using the current technique.

\subsection{Bounds on the Thurston-Bennequin number}\label{sec:TB}
One of the more useful application of the odd Khovanov homology is in
finding upper bounds on the Thurston-Bennequin number of Legendrian links.
Consider $\R^3$ equipped with the standard contact structure $dz-y\,dx$.
A link $K\subset\R^3$ is said to be {\em Legendrian} if it is everywhere
tangent to
the $2$-dimensional plane distribution defined as the kernel of this $1$-form.
Given a Legendrian link $K$, one defines its {\em Thurston--Bennequin number},
$tb(K)$, as the linking number of $K$ with its push-off $K'$ obtained using
a vector field that is tangent to the contact planes but orthogonal to the
tangent vector field of $K$. Roughly speaking, $tb(K)$ measures the framing of
the contact plane field around $K$. It is well-known that the TB-number can be
made arbitrarily small within the same class of topological links via
stabilization, but is bounded from above.

\begin{defin}\label{def:TB-bound}
For a given {\em topological} link $L$, let $\tbb(L)$, the {\em TB-bound}
of $L$, be the maximal possible TB-number among all the Legendrian
representatives of $L$. In other words, $\tbb(L)=\max_K\{tb(K)\}$, where
$K$ runs over all the Legendrian links in $\R^3$ that are topologically
isotopic to $L$.
\end{defin}

Finding TB-bounds for links has attracted considerable interest lately, since
they can be used to demonstrate that certain contact structures on $\R^3$ are
not isomorphic to the standard one. Such bounds can be obtained from the
Bennequin and Slice-Bennequin inequalities, degrees of HOMFLYPT and Kauffman
polynomials, Knot Floer homology, and so on (see~\cite{Ng-TB-bound}
for more details). The TB-bound coming from the Kauffman polynomial is usually
one of the strongest, since most of the others incorporate another invariant
of Legendrian links, the rotation number, into the inequality.
In~\cite{Ng-TB-bound}, Lenhard Ng used Khovanov homology to define a new
bound on the TB-number.

\begin{thm}[Ng~\cite{Ng-TB-bound}]\label{thm:Ng-TB}
Let $L$ be an oriented link. Then
\begin{equation}\label{eq:Ng-TB}
\tbb(L)\le\min\Bigl\{k\big|\bigoplus_{j-i=k}\CalH^{i,j}(L;R)\not=0\Bigr\}.
\end{equation}
Moreover, this bound is sharp for alternating links.
\end{thm}

This Khovanov bound on the TB-number is often better than those that were
known before. There are only two prime knots with up to 13 crossings for which
the Khovanov bound is worse than the one coming from the Kauffman
polynomial~\cite{Ng-TB-bound}. There are 45 such knots with at most 15
crossings.

\begin{figure}
\centerline{%
\vbox{\halign{\hfill#\hfill\cr\input{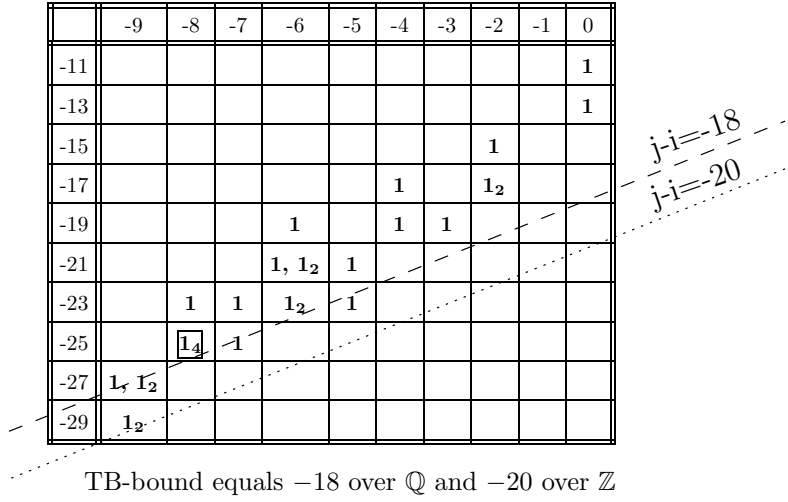}\cr
\noalign{\vskip -2\bigskipamount}
TB-bound equals $-18$ over $\Q$ and $-20$ over $\Z$\qquad\qquad\cr}}}
\caption{Khovanov TB-bound for the $(4,-5)$-torus knot}\label{fig:TB-bounds}
\end{figure}

\begin{example}
Figure~\ref{fig:TB-bounds} shows computations of the Khovanov TB-bound for the
$(4,-5)$-torus knot. The Khovanov homology groups in~\eqref{eq:Ng-TB} can be
used over any ring $R$, and this example shows that the bound coming from the
integral homology is sometimes better than the one from the rational one, due
to a strategically placed torsion. It is interesting to note that the integral
Khovanov bound of $-20$ is computed incorrectly in~\cite{Ng-TB-bound}. In
particular, this was one of the cases where Ng thought that the Kauffman
polynomial provides a better one. In fact, the TB-bound of $-20$ is sharp for
this knot.
\end{example}

The proof of Theorem~\ref{thm:Ng-TB} is based on the long exact
sequences~\eqref{eq:Kh-exact} and, hence, can be applied verbatim to the
reduced as well as odd Khovanov homology. By taking into account a grading
shift in the reduced homologies, we immediately get that
\begin{align}
\tbb(L)\le&-1+\min\Bigl\{k\big|\bigoplus_{j-i=k}\tCalH^{i,j}(L;R)\not=0\Bigr\}\\
\tbb(L)\le&-1+\min\Bigl\{k\big|
\bigoplus_{j-i=k}\tCalH_{\odd}^{i,j}(L;R)\not=0\Bigr\}.
\end{align}

As it turns out, the odd Khovanov TB-bound is often better than the even one.
In fact, computations with \KhoHo and Knotscape show that the odd Khovanov
homology provides the best upper bound on the TB-number among all currently
known ones for all prime knots with at most 15 crossings. More specifically,
the odd Khovanov TB-bound equals the Kauffman one on all the $45$ knots with
at most $15$ crossings where the latter is better than the even Khovanov
TB-bound. These knots are $12^n_{475}$, $13^n_{1708}$,
$\overline{14}^n_{9580}$, $\overline{14}^n_{9989}$, $14^n_{12208}$,
$14^n_{14433}$, $\overline{14}^n_{15458}$, $14^n_{18373}$, $14^n_{21980}$,
$14^n_{24190}$, $\overline{15}^n_{24518}$, $15^n_{34445}$, $15^n_{34827}$,
$\overline{15}^n_{37632}$, $15^n_{40088}$, $\overline{15}^n_{40854}$,
$\overline{15}^n_{42851}$, $\overline{15}^n_{49772}$,
$\overline{15}^n_{51188}$, $15^n_{51379}$, $15^n_{52894}$, $15^n_{52993}$,
$15^n_{53039}$, $15^n_{53226}$, $15^n_{54333}$, $15^n_{54502}$,
$15^n_{57673}$, $\overline{15}^n_{59594}$, $15^n_{62334}$, $15^n_{64226}$,
$15^n_{72658}$, $\overline{15}^n_{76240}$, $15^n_{79161}$, $15^n_{85661}$,
$\overline{15}^n_{92272}$, $\overline{15}^n_{100242}$,
$\overline{15}^n_{101483}$, $15^n_{103773}$, $15^n_{124839}$, $15^n_{124916}$,
$15^n_{132200}$, $15^n_{140988}$, $15^n_{144125}$, $\overline{15}^n_{145208}$,
and $15^n_{167608}$. On the other hand, there are many knots (several
thousands) for which the Khovanov TB-bound (both even and odd) is better than
the Kauffman one. The first such knot is $11^n_{20}$.

\begin{figure}
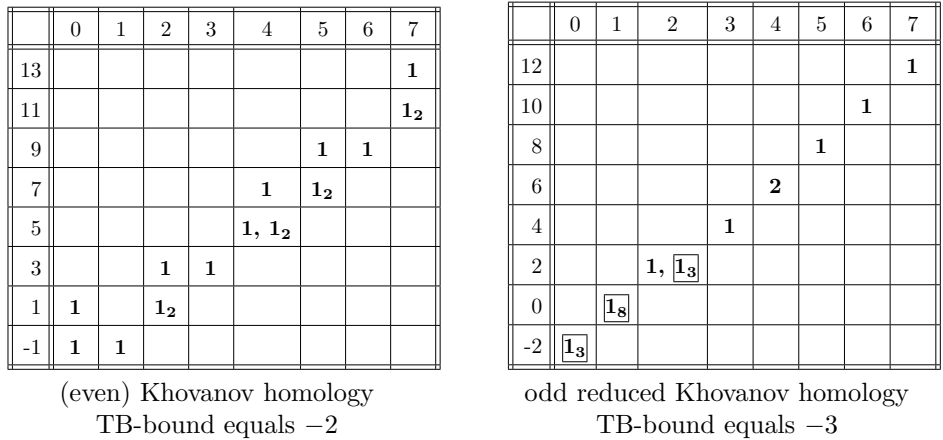

\centerline{%
\vbox{\halign{\hfill#\hfill\cr
\resizebox{0.47\hsize}{!}{\input{knot-12n_475-table}}\cr
(even) Khovanov homology\cr TB-bound equals $-2$\cr}}\qquad
\vbox{\halign{\hfill#\hfill\cr
\resizebox{0.47\hsize}{!}{\input{knot-12n_475-table-Odd}}\cr
odd reduced Khovanov homology\cr TB-bound equals $-3$\cr}}}
\caption{Khovanov TB-bounds for the knot $12^n_{475}$}\label{fig:12n_475-hom}
\end{figure}

\begin{example}
The odd Khovanov TB-bound is better than the even one and is equal to the
Kauffman one for the knot $12^n_{475}$, as shown in
Figure~\ref{fig:12n_475-hom}.
\end{example}

\subsection{Finding quasi-alternating knots}\label{sec:quasi-alt}
Quasi-alternating links were introduced by Ozsv\'ath and Szab\'o
in~\cite{OS-spectral} as a way to generalize the class of alternating links.

\begin{defin}\label{def:QA-links} The class $\CalQ$ of quasi-alternating links
is the smallest set of links such that
\begin{itemize}
\item the unknot belongs to $\CalQ$;
\item if a link $L$ has a planar diagram $D$ such that the two
resolutions of this diagram at one crossing represent two links,
$L_0$ and $L_1$, with the properties that $L_0,L_1\in\CalQ$ and 
$\det(L)=\det(L_0)+\det(L_1)$, then $L\in\CalQ$ as well.
\end{itemize}
\end{defin}

\begin{rem}
It is well-known that all non-split alternating links are quasi-alternating.
\end{rem}

The main motivation for studying quasi-alternating links is the fact that
the double branched covers of $S^3$ along such links are so-called $L$-spaces.
A $3$-manifold $M$ is called an {\em $L$-space} if the order of its first
homology group $H_1$ is finite and equals the rank of the Heegaard--Floer
homology of $M$ (see~\cite{OS-spectral}). Unfortunately, due to the recursive
style of Definition~\ref{def:QA-links}, it is often highly non-trivial to
prove that a given link is quasi-alternating. It is equally challenging to
show that it is not.

To determine that a link is not quasi-alternating, one usually employs the fact
that such links have homologically thin Khovanov homology over $\Z$ and Knot
Floer homology over $\Z_2$ (see~\cite{QA-links}).
Thus, $\Z$H-thick knots are not quasi-alternating. There are 12 such knots
with up to 10 crossings. Most of the others can be shown to be
quasi-alternating by various constructions. After the work of Champanerkar and
Kofman~\cite{Kofman-Co-QA}, there were only two knots left, $9_{46}$ and
$10_{140}$, for which it was not known whether they are quasi-alternating or
not. Both of them have homologically thin Khovanov and Knot Floer homology.

\begin{figure}
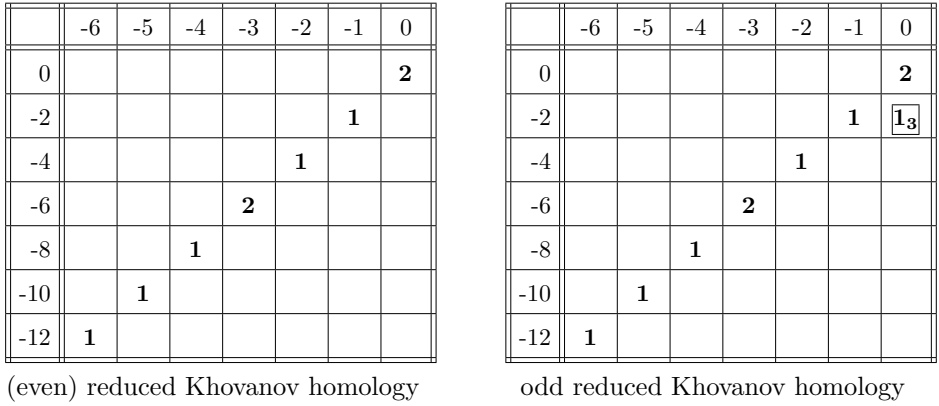

\centerline{%
\vbox{\halign{\hfill#\hfill\cr
\resizebox{0.47\hsize}{!}{\input{pretzel-3-3-n3-table-Red}}\cr
(even) reduced Khovanov homology\cr}}\qquad
\vbox{\halign{\hfill#\hfill\cr
\resizebox{0.47\hsize}{!}{\input{pretzel-3-3-n3-table-Odd}}\cr
odd reduced Khovanov homology\cr}}}
\caption{Khovanov homology of $9_{46}$, the $(3,3,-3)$-pretzel
knot}\label{fig:9_46-hom}
\end{figure}

\begin{figure}
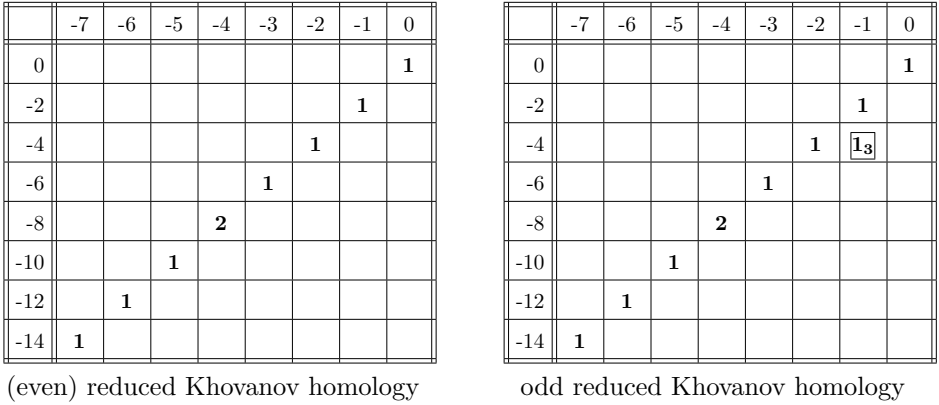

\centerline{%
\vbox{\halign{\hfill#\hfill\cr
\resizebox{0.47\hsize}{!}{\input{pretzel-3-4-n3-table-Red}}\cr
(even) reduced Khovanov homology\cr}}\qquad
\vbox{\halign{\hfill#\hfill\cr
\resizebox{0.47\hsize}{!}{\input{pretzel-3-4-n3-table-Odd}}\cr
odd reduced Khovanov homology\cr}}}
\caption{Khovanov homology of $10_{140}$, the $(3,4,-3)$-pretzel
knot}\label{fig:10_140-hom}
\end{figure}

As it turns out, odd Khovanov homology is much better at detecting
quasi-alternating knots. The proof of the fact that such knots are $\Z$H-thin
is based on the long exact sequences~\eqref{eq:Kh-exact} and, therefore, can
be applied verbatim to the odd homology as well~\cite{Khovanov-odd}.
Computations with \KhoHo show that the knots $9_{46}$ and $10_{140}$ have
homologically thick odd Khovanov homology and, hence, are not
quasi-alternating, see Figures~\ref{fig:9_46-hom} and~\ref{fig:10_140-hom}.

\begin{figure}
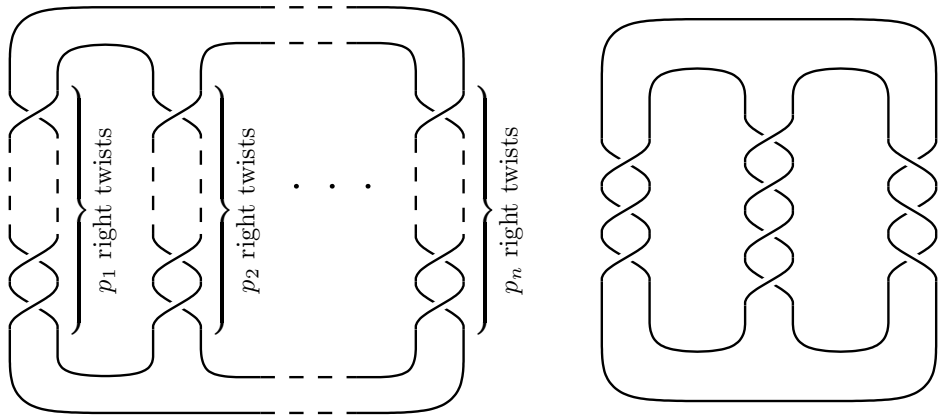

\centerline{$\vcenter{\hbox{\input{pretzel-link.pspdftex}}}\qquad\qquad
\vcenter{\hbox{\input{pretzel-3_4_n3.pspdftex}}}$}
\caption{$(p_1,p_2,\dots,p_n)$-pretzel link and $(3,4,-3)$-pretzel knot
$10_{140}$}\label{fig:pretzels}
\end{figure}

\begin{rem}\label{rem:odd-torsion}
It is worth mentioning that the knots $9_{46}$ and $10_{140}$ are $(3,3,-3)$-
and $(3,4,-3)$-pretzel knots, respectively (see Figure~\ref{fig:pretzels} for
the definition). Computations show that $(n,n,-n)$- and $(n,n+1,-n)$-pretzel
links for $n\le 6$ all have torsion of order $n$ outside of the main diagonal
that supports the free part of the homology. This suggest a certain
$n$-fold symmetry on the odd Khovanov chain complexes for these
pretzel links that cannot be explained by the construction.
\end{rem}

\begin{rem}
Joshua Greene has recently identified~\cite{Greene-QA} all quasi-alternating
pretzel links by considering $4$-manifolds that are bounded by the 
double branched covers of the links.
In particular, he found several knots that are not
quasi-alternating, yet are both H-thin and OH-thin. The smallest such knot
is $11^n_{50}$. 
\end{rem}

\subsection{Transversely non-simple knots}
Consider once again $\R^3$ equipped with the standard contact structure
$dz-y\,dx$. A knot $K\subset\R^3$ is said to be {\em transverse} if it is
everywhere transverse to the contact planes (cf. Section~\ref{sec:TB}).
Given a transverse knot $K$, one defines its {\em self-linking number},
$sl(K)$, as the linking number of $K$ with its push-off $K'$ along the vector
field $v=\frac{\partial}{\partial y}$. This vector field always lies in the
contact planes.

\begin{defin}\label{def:non-simple}
A topological knot $K$ is said to be {\em transversely non-simple}, if it has
two transverse representatives $K_1$ and $K_2$ such that $sl(K_1)=sl(K_2)$,
but $K_1\not=K_2$ as transverse knots.
\end{defin}

The first examples of transversely non-simple knots were discovered in~2003.
Birman and Menasco found a family of transversely non-simple
$3$-braids~\cite{BirmanMenasco-1,BirmanMenasco-2}, while Etnyre and Honda
proved that the $(2,3)$-cable of the trefoil is transversely non-simple as
well~\cite{Etnyre-Honda}. More recently, several families of transversely
non-simple knots were found using knot Floer
homology~\cite{Baldwin-transverse,Ng-Co-transverse,Vertesi-transverse}.

It turns out that many of these transversely non-simple knots have very
special odd Khovanov homology. Namely, they have only torsion and no rational
homology in the homological grading $0$. We call such knots {\em
zero-omitting}. There are 677 zero-omitting knots among all prime knots with
at most $15$ crossings. 10 of them have 12 crossings or less: $10_{132}$,
$11^n_{38}$, $12^n_{120}$,
$12^n_{199}$, $12^n_{200}$, $12^n_{260}$, $12^n_{475}$, $12^n_{523}$,
$12^n_{549}$, and $12^n_{673}$. Homology of one of these knots is shown in
Figure~\ref{fig:12n_475-hom}.
Seven out of $10$ of the zero-omitting knots with at most 12 crossings are
known to be transversely non-simple~\cite{Ng-private}. The three unknown cases
are $11^n_{38}$,  $12^n_{475}$, and $12^n_{673}$. On the other hand, the
$(2,3)$-cable of the trefoil, a transversely non-simple knot, is not
zero-omitting.


\begin{conjec}
All zero-omitting knots are transversely non-simple.
\end{conjec}

\end{document}